\documentclass[12pt]{amsart}
\newtheorem{theo}{THEOREM}[section] 
\newtheorem{lemma}[theo]{Lemma}

\newtheorem{prop}[theo]{Proposition}
\theoremstyle{remark}

\newcommand{\brref}[1]{(\ref{#1})}
\newcommand{\calo}{\mathcal {O}}
\newcommand{\tensor}{\otimes}

\newcommand{\Proj}[1]{\mathbb{ P}(#1)}

\newcommand{\restrict}[2]{{#1}_{\mid _{#2}}}
\newcommand{\oof}[2]{{\mathcal O}_{#1}({#2})}
\newcommand{\oofp}[2]{{\mathcal O}_{{\mathbb P}^{#1}}({#2})}

\newcommand{\sel}{(S, L)}

\newcommand{\Pin}[1]{{\mathbb P}^{#1}}
\newcommand{\taut}[1]{{\mathcal O}_{\mathbb{P}(#1)}(1)}
\newcommand{\tautof}[2]{{\mathcal O}_{\mathbb{P}(#1)}(#2)}

 \newcommand{\Sq}{{\mathcal S}_4}
 \newcommand{\Sqq}{{\mathcal S}^*_4}
 \newcommand{\St}{{\mathcal S}_3}
 \newcommand{\Sn}{{\mathcal S}_n}
 \newcommand{\SinSt}{\sel \in \St}
 \newcommand{\SinSq}{\sel \in \Sq}
 \newcommand{\SinSqq}{\sel \in \Sqq}
 \newcommand{\SinSn}{\sel \in \Sn}
 \newcommand{\psius}{\psi^*}

 \newcommand{\proj}[1]{\mathbb{P}(#1)}
\newcommand{\commtri}[6]{\begin{center}
\begin{picture}(8,6)
\put(1,5){$#1$} 
\put(7.3,5){$#2$}
\put(4.5,0.5){$#3$} 
\put(2,5){\vector(1,0){5}}
\put(2,4.5){\vector(2,-3){2}}
\put(7.5,4.5){\vector(-2,-3){2}}
\put(2.2,2.3){$#6$} 
\put(4.5,5.5){$#4$}
\put(6.6,2.5){$#5$} 
\end{picture}
\end{center}} 
 \title[On polarized  surfaces]{On polarized surfaces of low degree  whose adjoint bundles are not
spanned}
 \author[G. Besana ]{Gian Mario Besana}
 \address{Department of Mathematics- Eastern Michigan University-
 Ypsilanti MI 48197- U.S.A.}
 \email{ gbesana@@emunix.emich.edu }
 \author[S. Di Rocco] {Sandra Di Rocco}
\address{ Department of Mathematics- KTH - Royal Institute of Technology - S 100 44 Stockholm -
SWEDEN}
\email{sandra@@math.kth.se }
\begin{document}

 
 \begin{abstract}
 Smooth complex surfaces polarized with an ample and globally generated
 line bundle of degree three and four, such that the adjoint bundle is not globally
 generated, are considered. Scrolls of a vector bundle over a smooth curve are shown to be the
only examples in degree three. Two classes of examples in degree four are presented, one of which is shown
to characterize regular such pairs. A Reider-type theorem is obtained in which the assumption on the
degree of $L$ is removed.
 \end{abstract}
 
 \maketitle
 
 
 \section{introduction}
 
 Let $\sel$ be a  pair where $S$ is a smooth complex projective surface and $L$  is an ample line bundle on
$S.$ Let $K$ be the canonical bundle of $S.$ If $L$ is very ample, it is a classical result of Sommese and Van De
Ven's
\cite{so-v} that the adjoint linear system  $|K+L|$  is free, unless $\sel$ is a scroll  over a smooth curve (see
section
\ref{notation} for definitions) or one of the special surfaces with sectional genus zero. If $L$ is only ample
and spanned, Reider's theorem,
\cite{Re}, implies that $|K+L|$ is free unless   $\sel$ is a scroll, under the assumption that $L^2 \ge 
 5.$ This assumption on the degree of the polarization is essential for
 Reider's method, being equivalent to the condition for  Bogomolov's
 instability of suitable rank-two vector bundles on $S.$
 Understanding what happens {\it below Reider} is, in Fujita's words,
 \cite{fu} p.156, {\it an interesting but subtle problem}. A first attempt was made in \cite{BE2}, where the
case of  Kodaira dimension $Kod(S) \le 0$ for $L^2=3$ was treated.
 In this paper,  polarized pairs of degree $3$ and $4$ are considered where, again, $L$ is  ample and spanned
and $|K+L|$ is not free. Theorem \ref{s3empty} shows that scrolls are the only such surfaces when
$L^2=3.$   Only two examples of such pairs of  degree four, which are not scrolls, were known to the
authors. They could  be found for example in \cite{fu} and they are a Del Pezzo surface of degree $1,$
polarized with
$L=-2K,$ and an elliptic $\Pin{1}$-bundle $S=\proj{E}$ of invariant $e=-c_1(E)=-1,$ polarized by
$\tautof{E}{2}.$   Two large  families of examples in degree four are presented in sections \ref{thefamilyS0} 
and \ref{thefamilyI}, to which the above mentioned examples belong. The Del Pezzo example is generalized to
a family of double covers of quadric cones, while the elliptic $\Pin{1}$ bundle is generalized to a class of
quotients of products of hyperelliptic curves. The first  of these families is shown in Theorem
\ref{s4first} to characterize the regular such pairs. 
Combining the results presented in this work with \cite{so-v} and \cite{Re} the following theorem is
obtained:
\begin{theo}
Let $S$ be   a smooth projective complex surface and let $L$ be an ample and spanned line bundle on $S.$
$|K+L|$ is free if and only if $\sel$ is not one of the following pairs :
\begin{itemize}
\item[i)]  $ (\Pin{2}, \oofp{2}{1})$ or $(\Pin{2}, \oofp{2}{2});$ 
\item[ii)] A scroll over a smooth curve;
\item[iii)]  A double cover of a quadric cone in $\Pin{3},$ given by $|L|$,  ramified over the vertex and
the intersection of the cone with a  general surface of degree $2a+1;$
\item[iv)] $L^2=4,$ $q(S)>0,$ $h^0(L)=3,$  and for every
$x\in Bs|K+L|$ there exists  a unique
$C\in |L - x|$ such that
$C=A+B$ where $A$ and $B$ are irreducible, reduced, ample divisors, $AB=1$, $A \equiv B,$ $\oof{A}{B} =
\oof{A}{x},$ $h^0(A)=h^0(B)=1.$
\end{itemize}

\end{theo}
 Both authors are grateful to Andrew Sommese for his insight in the construction of the families of
examples of degree four, to Kristian Ranestad for helpful conversations and to the referee for suggested
improvements in the presentation of \brref{thefamilyS0}.  Both authors would like to thank the Max
Planck Institute f\"ur Mathematik in Bonn, Germany, and the Kungl Tekniska H\"ogskolan (Royal Institute
of Technology) in Stockholm, Sweden,  where this work was carried out.  The first author would like to
thank the Graduate School of Eastern Michigan University for its support through two  Graduate School
Research Support Awards.
 
 
 \section{Notation} \label{notation}
 
 Throughout this article $S$ denotes  a smooth, connected, projective surface defined over the complex field
${\mathbb C}.$ Its structure sheaf is denoted  by ${\mathcal O}_S$ and the canonical sheaf of
holomorphic $2$-forms on $S$ is denoted by $K_S$ or simply $K$ when the ambient surface is understood.
 For any coherent sheaf ${\mathcal F}$ on $S$,  $h^i({\mathcal F} )$ is the complex dimension
of $H^i(S,{\mathcal F}).$ 
Let $L$ be a  line bundle on $S.$ If $L$ is ample the pair $\sel$ is called a {\it polarized surface}.
The following notation and definitions  are used:\\
  $|L|$, the complete linear system associated with L;\\
 $Bs|L|$, the base locus of the linear system $|L|$;\\
$|L|$ is {\it free at a point} $x$ if $x \not \in Bs|L|.$
$|L|$ is {\it free} if $Bs|L|=\emptyset$ or equivalently if  $L$ is spanned, i.e.
generated by its global sections.\\
 $d = L^2,$ the degree of $L$;\\
$g=g(S, L)$, the {\it sectional genus} of $\sel$, defined by $2g-2=L (K_S+L).$ \\
$\Delta \sel = \Delta = 2+L^2-h^0(L),$ the Delta genus of $\sel$;\\ 
$\mathbf{F_e}$, the Hirzebruch surface of
invariant $e$ ;\\
A polarized surface $(S, L)$
is a {\em scroll} over a smooth curve $C$ if there exists a rank two vector bundle $E$ over $C,$ such that
$(S,L) = (\Proj{E},\taut{E}).$\\
$\sigma: \hat{S}=Bl_P S \to S,$ the blow up of a surface $S$ at a point $P.$\\
Cartier divisors, their associated line bundles and the invertible sheaves of their holomorphic sections are
used with no distinction. Mostly additive notation is used for their group. Given two divisors $L$ and $M$
we denote linear equivalence by $L \sim M$ and numerical equivalence by $L \equiv M.$
 When $S$ is a $\Pin{1}$-bundle over a curve with fundamental section $C_0$ and generic  fiber $f$ we
have
$Num(S) ={\mathbb Z}[C_0] \oplus {\mathbb Z}[f].$ \\ 
A standard argument, see for example \cite{bi-fa-la}, shows that a
polarized surface $\sel$ with $L$ ample and spanned is a scroll if  there exists an effective divisor $E\subset
S$ such that $E^2=0$ and $LE=1.$ This result will be used repeatedly in this work.
As in  \cite{BE2},  the   following notation will be used:
 
 \begin{xalignat}{2}
  \Sn =  \{ \sel |& L \text{ ample and spanned but not very ample,}\\ \notag
  &  L^2=n,\\ \notag
  &  |K+L| \text{ not free},\\ \notag
  &  \sel \text{ not a scroll}.\} 
 \end{xalignat}
 
  For $\SinSn,$ $ \psi$ will always denote the holomorphic map given by $|L|.$

%
%
 \section{The key tools} The assumption $L^2 \ge 5$ in Reider's theorem is equivalent to the
Bogomolov's instability of a suitable vector bundle, whose existence is guaranteed by
Cayley-Bacharach type conditions.  As such,  it is essential to his  method. On the other hand,  
Sommese's original argument in
\cite{SO1}, uses the  very ampleness of $L$ exclusively in order to satisfy the key requirement that
is highlighted in the following proposition. 
%
%

\begin{prop}
\label{key} Let $L$ be an ample and spanned line bundle on a smooth, projective, irregular surface
$S.$ Then
$|K+L|$ is free at
$x\in S$ if for any fixed tangent direction $\mathbf{v}$ at $x,$  there exists a curve $C \in |L - x|$ smooth at
$x,$ having $\mathbf{v}$ as tangent direction at $x$ and such that
$|\omega_C|$ is base point free. 
\end{prop}
\begin{proof} Let $\Lambda= |L- x|$ be the linear system of curves through $x.$ Let $C \in
\Lambda$ and let $\omega_C$ be its dualizing sheaf. Consider the following commutative diagram (for
details see Andreatta and Sommese, \cite{an-so2}).
\begin{xalignat}{2} 
 ..\to H^0(K+L) \stackrel{\alpha}{\rightarrow}H^0(C,\omega_C)& 
 \stackrel{\beta}{\rightarrow} H^1(S,K) \to 0 \label{diagAndrew} \\ 
  r \nwarrow & \ \   \nearrow \gamma \ \\ \notag  
  H^0(&S,\Omega^1_S)  
 \end{xalignat}
 The map $\gamma: H^0(S,\Omega^1_S) \to H^1(S,K),$ given by wedging with $c_1(L),$
 gives an isomorphism. Because $q\ge 1$
 it is $H^0(S,\Omega^1_S) \neq 0.$
 The above diagram gives 
 \begin{equation} 
 \label{directsumforKC}
 H^0(\omega_C) = Im(r) \oplus Im (\alpha)
 \end{equation} 
 If the cotangent bundle $\Omega^1_S$ is generated by its global
 sections at $x,$ then we can find two linearly independent holomorphic
 $1$-forms, $\eta_1$, $\eta_2,$ non vanishing at $x,$ and thus
 $\eta_1\wedge\eta_2$ would give an holomorphic two form non vanishing
 at $x$. Because $L$ is
 spanned, there exists  $\sigma \in H^0(L)$  such that $\sigma(x) \ne 0.$ Then
 $(\eta_1 \wedge \eta_2)\otimes \sigma$ is a section of $K+L$ which
 does not vanish at $x,$ and thus $K+L$ is spanned at $x.$
 We can then assume that the evaluation map $ev_x:H^0(\Omega^1_S)\to  H^0( \Omega^1_{S,x})$  is
not surjective. If $\dim(ev_x(H^0(\Omega^1_S))=0$ then every section of $\omega_C$ non 
 vanishing at $x$ is of the form $\alpha(s)$ with $s$ a section of
 $K+L$ non vanishing at $x$. Thus $K+L$ is spanned at $x$ since $|\omega_C|$ is base point free. Assume
finally that $\dim(ev_x(H^0(\Omega^1_S))= 1.$ Then there exists a tangent direction $\mathbf{v}$ such
that given $C \in |L|$  with tangent direction $\mathbf{v}$ at $x,$ for all $\omega \in H^0(S,
\Omega_S^1)$ it is $r(\omega)(x)=0.$ Base point freeness of $|\omega_C|$ and
\brref{directsumforKC} then give $|K+L|$ spanned at $x.$
\end{proof}

The above Proposition shows that it will be necessary to establish base point freeness of
$|\omega_C|$ for possibly singular curves on $S.$ This problem was studied for example by Catanese,
\cite{cat}. His results, together with results due to Francia, \cite{Fr},  have recently been
reinterpreted by Mendes Lopes,
\cite{margarida}, in a setting quite similar to ours.  For the convenience of the reader, two  results
from
\cite{margarida} are recalled:

%
%
 
 \begin{lemma} \cite[Theorem 3.1]{margarida}
\label{marga31} Let $C$ be a $1$-connected divisor on a smooth surface $S.$ Then a multiple point
$x \in C$ is a base point for $|K_S + C|$ if and only if $C$ decomposes as $C=A+B$ with:
\begin{itemize}
\item[a)] AB=1;
\item[b)] $x$ is a smooth point of $A$ and $\oof{A}{x} = \oof{A}{B}.$
\end{itemize} Furthermore if $x$ is a base point of $|K_S + C|$ then the decomposition above is such
that $A \cap B=\{x\}$ or $A \subset B.$
\end{lemma}
 
 \begin{lemma} \cite[Theorem 4.1]{margarida}
\label{marga41} Let $C$ be a $1$-connected divisor on a smooth surface $S.$ Let $x$ be a smooth
point of $C.$ Then $x$ is a base point of $|\omega_C|$ if and only if either $C \simeq \Pin{1}$ or $C$ is
reducible, the unique component $\Gamma$ to which $x$ belongs is a non singular rational curve 
and $C$ decomposes as $C= \Gamma + F_1 + \dots F_n,$ where the $F_i's$ are effective non-zero
divisors, satisfying:
\begin{itemize}
\item[i)] $F_i \Gamma=1$ for every $i;$
\item[ii)] $F_i F_j=0$ for $i \neq j;$
\item[iii)] $\oof{F_i}{F_k} \simeq \calo_{F_i}$ for $k < i.$
\end{itemize} Furthermore if $x$ is a base point of $|\omega_C|$ then $\Gamma$ is a fixed
component of
$|\omega_C|.$
\end{lemma}

%
%
The following simple facts will also be used:
 \begin{lemma}  \cite[Lemma 3.1]{BE2}
 \label{prel3}
 Let $\SinSt$ and let $C \in |L|$ be a smooth generic curve. 
 Then  
 \begin{itemize}
 \item[a)] $h^0(L)=3$, i.e $|L|$ expresses $S$ as a triple cover of $\Pin{2};$
 \item[b)]The restriction map $H^0(S,L) \to H^0(C,\restrict{L}{C})$ is onto, $q(S)\ge 1$ and $g=g \sel
\ge 2;$
 \item[c)] $h^0(K+L) > 0$ and $q < g;$
 \end{itemize}
 \end{lemma}
 \section{Polarized surfaces of degree three}
 \label{degree3}

%
%
 
 Following Mendes Lopes \cite{margarida} and others, see for example Catanese, \cite{cat}, the
investigation of the $n$-connectedness properties of curves in $|L|$ reveals crucial facts
 on the base locus of their dualizing linear system.
 \begin{lemma}
 \label{nconnected} Let $\SinSt$ and let $C\in |L|.$ Then $C$ is $2$-connected.
 \end{lemma}
 \begin{proof}
  Since $L$ is ample, $C$ is 1-connected. Assume there exist $C \in |L|$
 not $2$-connected. Then $C=A_1+A_2$ with $A_i$ effective and
 $A_1A_2=1.$ Since $LC=3$ and $L$ is ample, it must be $LA_i=1$ for one
 $i.$ Say $LA_1=1.$ Then $1=LA_1=(A_1+A_2) A_1= A_1^2+1$ and thus
 $A_1^2=0.$ But then $\sel$ is a scroll, contradiction. 
 \end{proof}
 %
 %

 
 \begin{lemma}
 \label{AB3}
 Let $\SinSt$  and let 
 $x \in Bs|K+L|.$ Let $\Lambda=|L- x|.$
 Then all the  $C \in \Lambda$ are  smooth at $x$
   and meet transversely at $x$.
 \end{lemma}
 \begin{proof}
 Let $C\in \Lambda.$ If $C$ were singular at $x$ then Lemma \ref{marga31} would
imply that 
$C$ is not
$2$-connected, which
 contradicts Lemma \ref{nconnected}. Therefore every $C \in \Lambda$ is smooth at $x.$
 The last part of the statement is a simple check in local coordinates,
 noticing that $\Lambda$ is a pencil.
 \end{proof}

 %
 %

 
 \begin{prop} 
 \label{wcfree}
 Let $\SinSt$ and let $C\in |L|.$ Then $|\omega_C|$ is base point free
 \end{prop}
 \begin{proof}
 By contradiction assume that $x\in C$ is a base point for
 $|\omega_C|$ and thus for $|K+L|.$ Lemma \ref{AB3} implies that $x$ is a smooth point  on $C.$
 Because $g(C)\ge 2,$ according to Lemma \ref{marga41},  $C$ must
 be reducible as $C= \Gamma + F_1 + \dots + F_n$ where $\Gamma$ is a smooth rational curve, $x
\in \Gamma,$ the $F_i$ are effective
 divisors, not necessarily irreducible,  such that $\Gamma F_i=1$ for all $i,$ and $F_i F_j = 0$
 for $i \ne j.$ Since $L$ is ample and
 $LC=3$ it must be $n\le 2.$ If $C= \Gamma + F_1$ the condition $\Gamma F_i=1$ violates the two
connectedness established in Lemma
 \ref{nconnected}.
 If $C= \Gamma + F_1 +F_2$ similarly  $F_1(\Gamma + F_2) = F_1 \Gamma = 1$ violates the
$2$-connectedness of
$C.$
 \end{proof}
 
 The following theorem is the central result of this section:

 %
 %


 \begin{theo}
\label{s3empty}
 $\St = \emptyset.$
 \end{theo}
 \begin{proof}
 By contradiction, assume there exists $\SinSt$ and let $x \in Bs|K+L|.$  Let $\Lambda= |L- x|$ be the pencil
of
 curves through $x.$ Let $C \in \Lambda$ and let $\omega_C$ be its
 dualizing sheaf. According to Lemma \ref{AB3} all $C \in \Lambda$ are smooth and transverse at
$x.$ Therefore it is possible to find a $C \in \Lambda$  having at $x$ any assigned tangent direction.
For such a  $C,$ the linear system $|\omega_C|$ is  base point free, according to Proposition \ref{wcfree} and
thus Proposition \ref{key} shows that $K+L$ is spanned at $x,$ contradiction.
\end{proof}


%
 %


 
 \section{Two classes of examples in $\Sq.$}
 In this section, the construction of two families of examples of  $\SinSq$ is presented.
 
 \subsection{A family $\mathcal{S}^0$ of regular polarized surfaces in $\Sq$}
 \label{thefamilyS0}
 Let $a\ge 1$ be an integer. Let  $Q$ be a rank $3$ quadric, in $\Pin{3},$ i.e. a cone with vertex $v$
over a smooth conic. Let  $\pi: {\mathbf F}_2 \to Q$ be the resolution of the vertex singularity, where $E$
denotes the exceptional divisor.  Let $B$ be the smooth  intersection of $Q$ with a general hypersurface of
degree $2a+1$ and ${\mathcal B}=\pi^{-1}(B) \cup E=\oof{{\mathbb{F}_2}}{(2a+2)E
 + (4a+2)f}=2{ \mathcal L}.$ Let  $\hat{\psi}$ be the  double cover of
 
 $\mathbb{F}_2$ given by ${\mathcal L}.$ Then we can consider the commutative diagram:
 \begin{center}
 \begin{tabular}{lcr}
 $S_a$& $\stackrel{\psi}{\rightarrow}$&$Q$\\
 $\sigma \uparrow$& & $\uparrow \pi$\\
 $\hat{S_a}$& $\stackrel{\hat{\psi}}{\rightarrow}$&${\mathbf F_2}$
 \end{tabular}
\end{center} where $\sigma:\hat{S_a} \to S_a$ is the contraction of the $(-1)$-curve ${\mathcal E}=
\hat{\psi}^{-1}(E)$ and $\psi: S_a\to Q$ is a double cover  with branch locus $B \cup {v}.$ Let
$L=\psi^*(\oof{Q}{1}).$ Notice that $L$ is ample and spanned, and $L^2=4.$  Fujita calls $(S_a,L)$ a
hyperelliptic manifold of type  $*II_a,$ \cite{fuhy}, and he shows that $Kod(S_1) =- \infty,$ while $Kod(S_a)
= 2$ for all $a \ge 2.$ The invariants of these surfaces are $q=0, g=2a,p_g=a(a-1)$ and $K^2=(2a-3)^2.$ 
Notice that $2 {\mathcal E} = \hat{\psi}^*(E).$ We have $$\sigma^*(L)=\hat{\psi}^*(E + 2f) =
2\mathcal{E}+ 2 \hat{\psi}^*(f).$$ As Fujita points out,
\cite{fuhy}  (4.6) (1), because  the line bundle
$\mathcal{E} + \hat{\psi}^*(f)$ on $\hat{S}_a$ is trivial on $\mathcal{E},$ there is a line bundle $H \in
Pic(S_a)$ such that $\sigma^*(H)= \mathcal{E} + \hat{\psi}^*(f)$ and therefore $L=2H,$ with $H$  ample.
Notice that $h^0(H)=2$ and thus the polarized pair $(S_a, H)$ is  a 
hypersurface of degree
$4a+2$ in the weighted projective space $\mathbb{P}(2a+1,2,1,1),$ according to \cite{fu}(6.21).
This easily implies that $K=(2a-3)H,$ which shows that  $S_1$ is a Del Pezzo surface of degree one with
$L=-2K,$ while
$K$ is ample and hence
$S_a$ minimal for $a \ge 2.$ 

 Let $x=\psi^{-1}(v).$ $K+L$ is spanned at $x$ if and only if
 $H^1(K+L - x)=0.$ This is equivalent to $H^1(\sigma^*(K+L) -
 {\mathcal E})= H^1(\hat{K} + \hat{L} - {\mathcal E})=0$ where $\hat{K}$ is the
 canonical bundle of $\hat{S}_a$ and $\hat{L}=\sigma^*(L)-{\mathcal E}={\mathcal E} +
2\hat{\psi}^*(f).$ By Serre duality 
 \begin{xalignat}{2}
   H^1(\hat{K} +\hat{L}
  - {\mathcal E})=&H^1(-\hat{L} + {\mathcal E})\\ \notag
  =&H^1(-2\hat{\psi}^*(f))\\ \notag
  =& H^1(-2f) \oplus H^1(-(2a+3)f - (a+1)E) =1 \ne 0.
 \end{xalignat}
 Therefore we obtained  a class of surfaces ${\mathcal S}^0 \subset \Sq.$
\subsection{ A family of irregular surfaces in $\Sq$}
\label{thefamilyI}
Let $C_1$ and $C_2$ be two copies of the same hyperelliptic curve of genus $q \ge 1.$ Let $X_q = C_1 \times C_2$
and let $\pi_i : X \to C_i$ be the projections onto the factors. Let $\iota : X_q \to X_q$ be the involution
$\iota (x,y)=(y,x)$ and let $S_q$ be the quotient $X_q/\iota.$ Let $p:X_q \to S_q$ be the resulting double
cover.

Consider a divisor $D=Q_1 + Q_2$ on $C_i,$ such that $|D|$ is a $g^1_2$ (the unique one if $q\ge 2$). Let
$\hat{L}_i=\pi^*_i(D)$ and let $\hat{L}=\hat{L}_1 + \hat{L}_2.$ It is $h^0(X_q,\hat{L})=4.$ Let
$H^0(X_q,\hat{L})^{\iota}$ be the subspace of global sections of $\hat{L}$  which are  $\iota$-invariant.
If $H^0(C_i,D)=<\sigma, \tau>$ then $H^0(X_q,\hat{L})^{\iota}=<\sigma \tensor \sigma, \tau 
\tensor \tau, \sigma \tensor \tau + \tau \tensor \sigma>.$ Now let $\hat{\mathcal{C}} \in |\hat{L}|$ and
consider the line bundle
$L$ on
$S_q$ associated to the divisor
$p(\hat{C}),$  so that $p^*(L)=\hat{L}.$  There is a natural isomorphism between global sections of $L$ and global
sections of
$\hat{L}$ which are $\iota$-invariant. Therefore $h^0(S_q, L)=3.$ From the construction it follows that $L$ is
ample and spanned, and $2L^2= \hat{L}^2 = (\hat{L}_1 + \hat{L}_2)^2 = 8$ so that $L^2=4.$

To show that $S_q \in \Sq$ it needs to be shown that $|K+L|$ is not free. Consider the smooth divisors
$L_i=p(\pi^*_1(Q_i) + \pi^*_2(Q_i)).$ It is $L_1 + L_2 \in |L|$ and $L_1 L_2 =1.$ If $Q_1\ne Q_2$
$L_1$ and $L_2$ meet transversely at a point $x=p(Q_1,Q_2).$  Therefore
$\restrict{(K+L)}{L_1}= K_{L_1}+\oof{L_1}{x}$  is not spanned at $x.$ Hence, $K+L$ is not spanned at
$x.$ Notice that by moving $Q_1+Q_2$ in the $g^1_2,$ one sees that $|K+L|$ has a fixed, rational
component. Observing that the branching locus of $p$ is isomorphic to $C_1$ and $C_2,$ one can
compute the invariants of $S_q:$
$$K^2=(q-1)(4q-9)\ \ \  \chi(\calo_{S_q})=\frac{(q-1)(q-2)}{2} \ \ \ h^1(\calo_{S_q})=q \ \ \ 
p_g(S_q)=\frac{(q-1)q}{2}.$$

Notice that when $q=1,$ $S_1$ is an elliptic $\mathbb{P}^1$-bundle.  If $P$ is a point on an elliptic curve
$C_1=C_2$ then taking $D=2P$ the above construction gives the same surface as the projectivization of 
the rank two vector bundle given by the extension $$0 \to \calo_{C_i} \to E \to \oof{C_i}{P} \to 0$$
where $L= \tautof{E}{2}.$

\section{Polarized surfaces of degree four}
%
%
It is quite natural to try to adopt the same approach as in 
section \ref{degree3} for $\SinSq.$ Unfortunately in this case $L$ is no
longer $2$-connected and this fact opens up an entirely different
scenario.

\begin{prop}
\label{prel4}
Let $\SinSq$ and let $C \in |L|$ be a smooth generic curve.
Then  $g(C)=g\sel \ge 2,$ $h^0(K+L) > 0,$ $q < g,$
and one of the   following occurs:
\begin{itemize}
\item[a)] $h^0(L)=4$ and  $\sel \in \mathcal{S}^0$
\item[b)] $h^0(L)=3$ and $|L|$ expresses $S$ as a quadruple cover of
$\Pin{2}.$
\end{itemize}
\end{prop}
\begin{proof}
Because $L$ is ample and spanned then $h^0(L) \ge 3$ and thus $\Delta
\le 3.$
Because it is always $\Delta \ge 0$ it is  $h^0(L) \le 6.$ If $\Delta =
0,1$
then $L$ must be very ample by \cite{fu} and thus $\sel \not \in \Sq$ by
\cite{so-v}. Therefore we have $\Delta\ge 2.$ If $g=0$ then $\Delta=0$ and thus $\sel \not \in \Sq.$ If $g=1$,
 being
$\Delta
\ge 2$ then $\sel$ should be an elliptic scroll, contradiction. Therefore $g\ge 2.$ 
As in \cite{BE2}, Lemma 3.1, it follows that $q<g$ and $h^0(K+L)>0.$ 

Assume $\Delta =2$, i.e. $h^0(L)=4.$
Then
$\psi : S
\to
\Pin{3}$
and $4=L^2=\deg \psi(S)\, \deg \psi.$
If $\deg \psi(S)=1$ then $\psi(S)=\Pin{2}$ and thus $h^0(L)=3$,
contradiction.
If $\deg \psi(S) =4$ then $\psi$ is a birational map onto a  quartic
hypersurface.
According to Fujita,
\cite{fu}, either $(S, L)$ is the blow up $\pi: S \to  \Pin{1} \times
\Pin{1}$ at 8
points polarized by $L=\pi^*(\oof{  \Pin{1} \times \Pin{1}}{2,3})-E_1
\dots
-E_8$ in which case $K+L= \pi^*(0,1)$ which is spanned, contradiction, or
$L$
is very ample which is also a contradiction.
If $\deg \psi(S) = 2$ then $\psi$ is a double cover of a quadric surface $Q$ in $\Pin{3}.$
Since $Q$ must be irreducible, either $Q$ is smooth or  it is
a cone over a smooth conic.
Assume $Q$ is smooth. Then $L= \psius (\oof{Q}{1,1})$ and
$K=\psius(\oof{Q}{-2+\frac{a}{2},-2+\frac{b}{2}})$ where
either $a,b$ are even, positive
integers, i.e.
$a
\ge 2, b\ge 2$,or $a=0, b\ge 2,$ since the branch locus of
$\psi$ must be effective and $Pic(Q)$ has no $2$-torsion.
Therefore $K+L= \psius(\oof{Q}{\frac{a-2}{2},\frac{b-2}{2}})$ and
$h^0(K+L)= h^0(Q,\oof{Q}{\frac{a-2}{2},\frac{b-2}{2}})
+h^0(Q,\oof{Q}{-1,-1})=h^0(Q,\oof{Q}{\frac{a-2}{2},\frac{b-2}{2}}).$
If $a$ and $b$ are both positive, $K+L$ is spanned, contradiction. If
$a=0$ then $h^0(K+L)=0$ which contradicts $h^0(K+L)>0$ established in the  above
paragraph, so in this case $\sel \not \in \Sq.$

Now assume $Q$ is a rank $3$ quadric, i.e. a cone with vertex $v$ over a
smooth
conic.  It follows from \cite{fuhy}, section 4,  that
these surfaces exist and they belong to the family $\mathcal{S}^0,$ described in \brref{thefamilyS0}

We are in case $b)$ if $\Delta=3$, i.e $h^0(L) = 3.$

\end{proof}

%
%
In view of Proposition \ref{prel4} the following notation will be
used in the sequel:
$$\Sqq:=\{\sel \in \Sq | h^0(L)=3.\}$$

\begin{lemma}
\label{AB}
Let  $\SinSqq$  and let
$x \in Bs|K+L|.$ Let $\Lambda=|L- x|.$
\begin{itemize}
\item[a)] If $C \in \Lambda$ is singular at $x$ then  $C=A+B$ where $A$ and $B$ are effective, ample,
irreducible, reduced, divisors with $A\equiv B$, $A^2=B^2=1,$ $LA=LB=2,$
$L\equiv 2A
\equiv  2B,$ $h^0(A)=h^0(B)=1,$ $g(A)=g(B)=\frac{g}{2}.$ Either $A$ and $B$ meet
transversely at $x$ or
$A=B$ and $C$ was not reduced;
\item[b)] There
exists a smooth  $C \in \Lambda.$
\end{itemize}
\end{lemma}
\begin{proof}
Since $L$ is ample, it is $1$-connected, thus \ref{marga31} implies
that a curve $C \in \Lambda$ singular at $x$ must be of the form $C=A+B$ where
$A$ and $B$ are effective and $AB=1.$  It is then
$4=L^2=(A+B)^2 = A^2 + B^2 + 2$ so that $A^2+B^2 = 2.$ Assume $A^2\le
0$, then
$B^2 \ge 2.$  Then the Hodge Index Theorem applied to $L$ and $B$
gives $(LB)^2
\ge 4B^2\ge 8$ i.e. $LB \ge 3.$  Since $L$ is ample and $L(A+B)=4$ it
must be $LB \le
3$ and so it must be $LB=3$ and thus $LA=1.$ Then $1=LA=(A+B)A=A^2+1$
gives
$A^2=0.$ Thus  $\sel$ must be a scroll, contradiction.
 Therefore it must be $A^2 > 0$ and similarly $B^2 >0,$
which
means that it must be $A^2=B^2=1,$ $LA=LB=2.$  The Hodge Index Theorem
now
gives $L \equiv 2A \equiv 2B$ and thus both $A$ and $B$ are ample
divisors, numerically equivalent to each other.
Since $AB=1$ it also follows that $A$ and $B$ are both irreducible and
reduced.

Because $L$ is ample and spanned, the image of the restriction map
$H^0(S,L)\to H^0(A,\restrict{L}{A})$ must be at least
two-dimensional. The sequence $0\to B \to L \to \restrict{L}{A},$
recalling that $h^0(L)=3,$
shows that $h^0(B)=1.$ Exchanging the role of $A$ and $B$ in the
argument we similarly get $h^0(A)=1.$
Computing the genus of $A$ and $B$ we have $$2g(A)-2=
(K+A)A=KA+1=\frac{KL}{2}+1 = g-2.$$
Therefore $g(A)=g(B) =\frac{g}{2}$.The last statement of part b)
follows from \ref{marga31}.

 To prove part $b)$  notice that $\Lambda$ is  a pencil, because $L$ is spanned.
If there is a smooth
 $C \in \Lambda$, part a) is proven, so assume that all $C \in \Lambda$
 are singular somewhere. Bertini's theorem
 implies that all curves in a Zariski open subset $W$ of $\Pin{1}$
  are smooth away from  $Bs(\Lambda).$
 Assume there exists a $C_1 \in W,$ not singular at $x$ and
 singular at $y_1 \in Bs (\Lambda),$ where obviously $y_1 \ne x.$
 If all other curves in $W$ are singular at $x$ then any two of them would have
 intersection $C_j C_k = L^2 \ge 5,$ contradiction. So assume $C_2\in
 W$ is smooth at $x.$ But then either $C_2$ is singular at $y_1$
 or at some other base point of $\Lambda$ and again $C_1 C_2 = L^2 \ge
 5,$ contradiction. Therefore all curves in $\Lambda$ are singular at
 $x,$ $Bs(\Lambda)=\{x\},$  and all curves in $\Lambda$ are reducible according to part a). Let
$\sigma: \hat{S}=Bl_x(S)\to S$ be the blow up of $S$ at $x.$ There is a map $\phi:\hat{S} \to \Pin{1}$
that, according to part $a),$ must factor through a double cover $\alpha$:
\commtri{\hat{S}}{\mathfrak{C}}{\Pin{1}}{\beta}{\alpha}{\phi}
Because the exceptional divisor of $\sigma$ dominates $\mathfrak{C},$ it must be $\mathfrak{C}=\Pin{1},$
but this contradicts $h^0(A)=h^0(B)=1.$
\end{proof}

The above Lemma allows us to give a detailed local picture of the linear system of curves in $|L|$ passing
through a base point of the adjoint system.
\begin{prop}
\label{localpic}
Let $\sel \in \Sq^*,$ let $x \in Bs|K+L|$ and let $\Lambda=|L- x|.$ Then there exists
exactly one $C \in \Lambda$ which is singular at $x$ and thus reducible as in Lemma \ref{AB} a), while all
other curves in $\Lambda$ are smooth at $x$ and have there the  same tangent.
\end{prop}
\begin{proof}
According to Lemma \ref{AB} there exists at least one smooth $C_1\in\Lambda.$ Assume there exists
another $C_2 \in \Lambda$ which is smooth at $x$  and meets $C_1$ transversely at $x.$ The existence of
these  two curves implies that every curve in $\Lambda$ is smooth at $x$ and allows us to find  a curve in
$\Lambda,$ smooth at $x,$ with any given tangent direction. Therefore Proposition \ref{key}
would give a contradiction if $|\omega_C|$ were free at $x$ for all $C \in \Lambda.$ Therefore, according to
Lemma \ref{marga41} and because $g(L)\ne 0,$ there exists a curve $C\in \Lambda$ such that $C=
\Gamma + F_1 + \dots F_n$ where $\Gamma$ is a nonsingular rational curve, passing through $x,$ which
is a fixed component for $|\omega_C|.$ Then $\Gamma$ is also a fixed component for $|K+L|.$ Lemma
\ref{AB}  implies that $n=1,$ $C= \Gamma+F_1,$ $F_1$ is irreducible and $g(F_1)=g(\Gamma)=0$ which
implies $g(L)=0,$ which is a contradiction.
Therefore every other  curve in $\Lambda$ which is smooth at $x$ must have the same tangent as $C_1.$
It is now a simple check in local coordinates to see that a pencil of curves through a point, that contains
smooth curves having the same tangent at the base point, contains only one singular element.
\end{proof}

We conclude this section  by showing that the family $\mathcal{S}^0$ characterizes the regular pairs in
$\Sq.$ We add a simple result on the relative minimality of pairs in $\Sq.$
%
%
\begin{theo}
\label{s4first}
Let $\SinSq.$ Then either $S \in \mathcal{S}^0$ or $q(S)\ge 1$ and then the map
given by $|L|$ expresses $S$ as a quadruple cover of $\Pin{2}.$
\end{theo}
\begin{proof}

Because of  Proposition \ref{prel4} to complete
the
proof  it must be shown that for every $\SinSqq,$ it is $q >0.$ Assume $\SinSqq$ and $q=0.$ Let
$x \in Bs
|K+L|.$ Since
$L$ is spanned, there is  a pencil $\Lambda$ of curves $ C \in |L|$
passing through $x.$ Lemma \ref{AB} guarantees that 
there exists  a smooth $C \in \Lambda.$ The fact that $g\ge 2$ and
the surjectivity of
$H^0(K+L) \to H^0(K_C)$ give $K+L$ spanned at $x$, contradiction. 
\end{proof}
%
%
\begin{lemma}
\label{relmin}
Let $\SinSq.$ Then $\sel$ is relatively minimal, i.e. there does not
exist any
$(-1)$-curve $E \subset S$ such that  $LE=1.$
\end{lemma}
\begin{proof}
Assume there exists a $(-1)$-curve $E$ such that $LE=1.$ Let $\sigma : S
\to S'$
be the contraction of $E.$ Let $L'$ be an ample line bundle on $S'$ such
that
$L=\sigma^*(L') - E.$ Then $(L')^2=5$ and therefore Reider's theorem
gives $K'+L'$
spanned (see the analogue argument in Lemma 3.7 of \cite{BE2}) and
thus
$K+L=\sigma^*(K'+L')$ also spanned, contradiction.
\end{proof}

%
%
\section{Proof of the Reider-type Theorem 1.1}
 Combining the material  presented above with Sommese-Van De Ven and Reider's previous results we
can now prove the  theorem presented in the introduction.
\begin{proof}
If  $L^2\ge 5$ the statement is due to Reider, \cite{Re}.
If $L^2=3,4$ the statement collects the above  Theorem \ref{degree3}, Theorem \ref{s4first},
Proposition \ref{localpic}.
If $L^2=1$ then $L,$ being ample and spanned,  is very ample and the result is due to Sommese and Van de
Ven, \cite{so-v}.
Let  $L^2=2.$ Because $L$ is ample and spanned, and $\Delta=4-h^0(L) \ge 0,$ it must be $h^0(L)=3,4.$  If
$h^0(L)=4$ then
\cite{fu} gives
$L$ very ample and the statement is due to Sommese and Van De Ven. If $h^0(L)=3$ then 
$|L|$ gives a double cover
$\psi:S\to \Pin{2}.$ Let
$\oofp{2}{a}$ be the line bundle such that the branching locus of the cover is a curve in $|\oofp{2}{2a}|.$
Then $K+L=
\psi^*(\oofp{2}{a-2}).$ Therefore either $K+L$ has no sections or it is spanned, thus $\sel \not \in \Sq.$
\end{proof}
 
 %
 
 
 %


 \end{document}